\documentclass[12pt,a4paper]{amsart}
\usepackage{amsmath,amssymb,amscd}
\input xy
\xyoption{all}

\newcommand{\ra}{\rightarrow}

\newcommand{\CC}{\mathbb C}
\newcommand{\KK}{\mathbb K}
\newcommand{\ZZ}{\mathbb Z}

\newcommand{\PP}{\mathbb P}

\newcommand{\cO}{\mathcal{O}}
\newcommand{\cM}{\mathcal{M}}
\newcommand{\cN}{\mathcal{N}}

\theoremstyle{plain}
\newtheorem{theorem}{Theorem}[section]
\newtheorem{lemma}[theorem]{Lemma}
\newtheorem{proposition}[theorem]{Proposition}
\newtheorem{cor}[theorem]{Corollary}
\newtheorem{remark}[theorem]{Remark}

\begin{document}
\title[Genus Zero III]{Coherent systems of Genus 0\\ 
III: Computation of flips for $k=1$}

\author{H. Lange}
\author{P. E. Newstead}

\address{H. Lange\\Mathematisches Institut\\
              Universit\"at Erlangen-N\"urnberg\\
              Bismarckstra\ss e $1\frac{ 1}{2}$\\
              D-$91054$ Erlangen\\
              Germany}
              \email{lange@mi.uni-erlangen.de}
\address{P.E. Newstead\\Department of Mathematical Sciences\\
              University of Liverpool\\
              Peach Street, Liverpool L69 7ZL, UK}
\email{newstead@liv.ac.uk}
\thanks{Both authors are members of the research group VBAC (Vector Bundles on Algebraic Curves). They were   
         supported  by the Forschungsschwerpunkt ``Globale Methoden in der komplexen Analysis'' of the DFG. 
         The second author would like to thank the Mathematisches Institut der Universit\"at 
         Erlangen-N\"urnberg for its hospitality}

\begin{abstract}
In this paper we continue the investigation of coherent systems of type $(n,d,k)$ on the projective line which are stable with respect to
some value of a parameter $\alpha$. We consider the case $k=1$ and study the variation of the moduli spaces with $\alpha$.
We determine inductively the first and last moduli spaces and the flip loci, and give an explicit description 
for ranks 2 and 3. We also determine the Hodge polynomials explicitly for ranks 2 and 3 and in certain cases for arbitrary rank.

\end{abstract}
\maketitle

\section{Introduction}

A {\it coherent system of type $(n,d,k)$} on a smooth projective curve $C$ over an algebraically closed field is by 
definition a pair $(E,V)$ with $E$
a vector bundle of rank $n$ and degree $d$ over $C$ and $V \subset H^0(E)$ a vector subspace of dimension $k$. 
For any real number $\alpha$, the {\it $\alpha$-slope} of a coherent system $(E,V)$ of type $(n,d,k)$ is defined by
$$
\mu_{\alpha}(E,V) := \frac{d}{n} + \alpha \frac{k}{n}.
$$
A {\it coherent subsystem} of $(E,V)$ is a coherent system $(F,W)$ such that $F$ is a subbundle of $E$ and $W \subset V \cap H^0(F)$.
A coherent system $(E,V)$ is called 
{\it $\alpha$-stable} ($\alpha$-{\it semistable}) if
$$
\mu_{\alpha}(F,W) < \mu_{\alpha}(E,V) \ \ (\mu_{\alpha}(F,W) \le \mu_{\alpha}(E,V))
$$
for every proper coherent subsystem $(F,W)$ of $(E,V)$. For every $\alpha$ there exists a moduli space $G(\alpha;n,d,k)$ of 
$\alpha$-stable coherent systems of type $(n,d,k)$.

In two previous papers \cite{ln} and \cite{ln1}, we computed in particular the precise conditions for existence of $\alpha$-stable 
coherent systems of type $(n,d,k)$ for $k=1,2,3$ on a curve of genus 0. In this paper we consider the case $k=1$ and investigate 
the relationship between the moduli spaces $G(\alpha;n,d,1)$ as $\alpha$ varies. As in \cite{ln1} our methods depend on 
the study of the ``flips'' which occur at ``critical values'' of $\alpha$.

There are only finitely many critical values of $\alpha$. Between any 2 consecutive critical values the moduli spaces 
$G(\alpha;n,d,1)$ do not change. To get from one moduli space to the next one crossing a critical value we have to delete a 
certain closed subvariety and insert another one in its place. Following \cite{bgn} we call this process a {\it flip}. The closed 
subvarieties are called {\it flip loci}.

After describing the general set up in section 2, we give an inductive determination of the flip loci in section 3.
This leads to a description of all moduli spaces as disjoint unions of locally closed subvarieties determined by 
the first moduli space and the flip loci (see Theorem \ref{thm3.5} for details).
In section 4 we determine inductively the first and last moduli spaces. 
In section 5 we introduce the Hodge polynomials and show how they can be explicitly computed for our moduli spaces.
In particular, the Hodge numbers $h^{p,q}$ are always zero for $p \neq q$.
In rank 2 our flips coincide precisely with those of Thaddeus \cite{th} and we give a complete analysis of this in section 6,
including explicit formulae for the Hodge polynomials of all moduli spaces.
In section 7 we consider coherent systems for rank $n \geq 3$ and determine the Hodge polynomials of the first moduli space 
for $t \leq 2$, where $t$ is defined by $d = na -t, 0 \leq t \leq n-1$. In particular, for $n=3$, this covers all cases 
and we show explicitly how each flip affects the Hodge polynomial.

Our original intention was to compute the Poincar\'e polynomials. The idea of using the Hodge polynomials which have 
stronger additive properties arose from the second author's reading of \cite{mov}. We are grateful to the authors. 

We work throughout on the projective line $\PP^1$ defined over an algebraically closed field $\KK$. For the discussion of 
Hodge polynomials $\KK = \CC$.

\section{The set up}

We consider the moduli spaces $G(\alpha;n,d,1)$ of $\alpha$-stable coherent systems on $\PP^1$ of type $(n,d,1)$.
When $n=1$, we can describe these spaces completely.

\begin{proposition} \label{prop1.1}
The moduli space $G(\alpha;1,d,1)$ is independent of $\alpha > 0$ and is isomorphic to $\PP^d$ when $d \geq 0$ 
and empty for $d < 0$.
\end{proposition}
 
\begin{proof}
For $n=1$ the $\alpha$-stability condition is just $\alpha > 0$. Moreover $h^0(\cO(d)) =0$ if $d < 0$ and  
$= d+1$ if $d \geq 0$. The result follows.
\end{proof} 

In view of this proposition we usually write $G(1,d,1)$ in place of $G(\alpha;1,d,1)$.

Suppose now that $n \geq 2$. Then we write as in \cite{ln} (and putting $k=1$),
\begin{equation} \label{eqn0}
d=na-t \quad \mbox{and} \quad a = l(n-1) + m + t.
\end{equation}
with $0 \leq t < n,\; l \geq 1$ and $0 \leq m < n-1$.
\begin{lemma} \label{lem1}
$d = (n-1)(ln+t) +mn$ and this determines $ln+t$ and $m$ subject to the condition $0 \leq m \leq n-2$.
\end{lemma}

\begin{proof}
\begin{eqnarray*}
d= na-t & = & na -nt +(n-1)t \\
& = & n(l(n-1) + m) + (n-1)t \\
& = & (n-1)(ln+t) + mn.
\end{eqnarray*}
\end{proof}

According to \cite[Theorem 5.1]{ln} there exist $\alpha$-stable coherent systems of type $(n,d,1)$ 
precisely when $\alpha$ is in the range
$$
t < \alpha < \frac{d}{n-1} - \frac{mn}{n-1}.
$$
This is equivalent to
\begin{equation} \label{e1}
t < \alpha < ln + t.
\end{equation}
For any allowable critical data set $(\alpha_c,n_1,d_1,k_1,n_2,d_2,k_2)$ defined as in \cite{ln1}, we have by assumption that
$$
\frac{k_1}{n_1} > \frac{k_2}{n_2}.
$$
In our case, this means
$$
k_1 = 1, \quad k_2=0
$$
and hence \cite[equation (5.2)]{ln1}
$$
n_2=1.
$$ 
Now write 
$$
d_2=a+e.
$$
By \cite[equation (2.9)]{ln1} it follows that
$$
\alpha_c = \alpha_e := en + t \quad \mbox{with} \quad 1 \leq e \leq l-1.
$$
This implies that for every critical value there is exactly one allowable critical data set.
In particular, if $l=1$, there are no critical values.

For any integer $e$ with $0 \leq e \leq l-1$ we write for brevity
$$
G_e := G(\alpha;n,d,1)
$$
for any $\alpha$ in the range $en+t < \alpha < (e+1)n+t$. 
If $(E,V) \in G_e$, we say also that $(E,V)$ is {\it $\alpha_e^+$-stable}. Similarly for $(E,V) \in G_{e-1}$ 
we say $(E,V)$ is {\it $\alpha_e^-$-stable}.

We denote by $G_e^+$ the flip locus consisting of those coherent
systems $(E,\KK) \in G_{e}$ which do not belong to $G_{e-1}$. 
As we shall see, these are given by exact sequences  
\begin{equation} \label{eq11}
0 \ra (\cO(a+e)^r,0) \ra (E,\KK) \ra (E_r,\KK) \ra 0
\end{equation}
for which no direct factor $(\cO(a+e),0)$ splits off and $(E_r,\KK)$ does not contain a subsystem $(\cO(a+e),0)$.

For a fixed $(E_r,\KK)$ with this property these extensions are described by the Grassmannian $Gr(r,Ext^1((E_r,\KK),(\cO(a+e),0)))$,
where according to \cite[equations (2.6) and (3.3)]{ln1},
$$
\dim \mbox{Ext}^1((E_r,\KK),(\cO(a+e),0)) = C_{12}
$$
with
\begin{eqnarray}  \label{e3}
C_{12} & = & -(n-r) - (a+e)(n-r) + (d-ra -re) + a + e + 1 \nonumber\\
& = & a-t - (n-1)(e+1) + r\\
& = & (l-e-1)(n-1) + m + r \geq r. \nonumber
\end{eqnarray} 

Similarly, for $1 \leq e \leq l$, we denote by $G_e^-$ the flip locus consisting of those coherent
systems $(E,\KK) \in G_{e-1}$ which do not belong to $G_e$. As we shall see, these are given by exact sequences 
\begin{equation} \label{eq1}
0 \ra (E_r,\KK) \ra (E,\KK) \ra (\cO(a+e)^r,0) \ra 0
\end{equation}
for which no direct factor $(\cO(a+e),0)$ splits off and $(E_r,\KK)$ does not admit any quotient system $(\cO(a+e),0)$.
For a fixed $(E_r,\KK)$ with this property these extensions are described by the Grassmannian $Gr(r,Ext^1((\cO(a+e),0),(E_r,\KK)))$,
where according to \cite[equations (2.7) and (3.4)]{ln1},
$$
\dim \mbox{Ext}^1((\cO(a+e),0),(E_r,\KK)) = C_{21}
$$
with
\begin{eqnarray} \label{e5}
C_{21} & = & -(n-r) + (a+e)(n-r) - (d-ra -re)  \nonumber\\
& = & -(n-r) + en + t\\
& \geq &  r. \nonumber
\end{eqnarray}  

Note that
\begin{equation}  \label{eqn7}
G_{e-1} \setminus G_e^- = G_e \setminus G_e^+.
\end{equation}

\section{Inductive determination of flip loci}

In this section we describe an inductive procedure to determine $G_e^-$ and $G_e^+$ for any $1 \leq e \leq l-1$. 
For this we need the following lemma

\begin{lemma} \label{lem2.1} Suppose $\alpha > 0$ and $b \in \ZZ$.

{\em (a)} Consider extensions of the form
\begin{equation} \label{eq2.1}
0 \ra (\cO(b)^r,0) \ra (E,\KK) \ra (E_r,\KK) \ra 0
\end{equation}
with $\mu_{\alpha}(E_r,\KK) = b$ and $(E_r,\KK)$ $\alpha$-semistable. Then $(E,\KK)$ is $\alpha^+$-stable 
if and only if 

{\em (i)} $(E_r,\KK)$ is $\alpha^+$-stable,

{\em (ii)} no direct factor $(\cO(b),0)$ splits off {\em (\ref{eq2.1})}.\\
If the class of {\em (\ref{eq2.1})} is given by $(e_1, \ldots, e_r)$ with $e_i \in {\rm Ext}^1((E_r,\KK),(\cO(b),0))$,
then {\em (ii)} is equivalent to 

{\em (ii$'$)} $e_1, \ldots, e_r$ are linearly independent.\\

{\em (b)} Consider extensions of the form
\begin{equation} \label{e2}
0 \ra (E_r,\KK) \ra (E,\KK) \ra (\cO(b)^r,0) \ra 0
\end{equation}
with $\mu_{\alpha}(E_r,\KK) = b$ and $(E_r,\KK)$ $\alpha$-semistable. Then $(E,\KK)$ is $\alpha^-$-stable 
if and only if 

{\em (i)} $(E_r,\KK)$ is $\alpha^-$-stable,

{\em (ii)} no direct factor $(\cO(b),0)$ splits off {\em (\ref{e2})}.\\
If the class of {\em (\ref{e2})} is given by $(e_1, \ldots, e_r)$ with $e_i \in {\rm Ext}^1((\cO(b),0),(E_r,\KK))$,
then {\em (ii)} is equivalent to 

{\em (ii$'$)} $e_1, \ldots, e_r$ are linearly independent.
\end{lemma}

\begin{proof}
(a): Suppose first $(E,\KK)$ is $\alpha^+$-stable. Then (ii) is obvious. For the proof of (i) suppose 
$(F,W) \subset (E_r,\KK)$ contradicts $\alpha^+$-stability. Then it also contradicts $\alpha$-stability and so
$\mu_{\alpha}(F,W) = \mu_{\alpha}(E_r,\KK)$. If $W=0$, the function $\mu_{\alpha}(E_r,\KK) - \mu_{\alpha}(F,0)$ is strictly 
increasing with $\alpha$. So $\mu_{\alpha^+}(E_r,\KK) > \mu_{\alpha^+}(F,0)$, a contradiction. Hence we must have $W = \KK$.
Then the pullback $(G,\KK)$ of $(F,\KK)$ to $(E,\KK)$ has $\mu_{\alpha}(G,\KK) = \mu_{\alpha}(E_r,\KK)$. This implies 
$\mu_{\alpha^+}(G,\KK) > \mu_{\alpha^+}(E,\KK)$, a contradiction.

Suppose, conversely, (i) and (ii) hold and $(G,W)$ is a subsystem of $(E,\KK)$ with 
$\mu_{\alpha^+}(G,W) \geq \mu_{\alpha^+}(E,\KK)$. Then we must have $\mu_{\alpha}(G,W) = \mu_{\alpha}(E,\KK)$, because 
$(E,\KK)$ is certainly $\alpha$-semistable. As above this implies $W = \KK$. The image of $(G,\KK)$ in $(E_r,\KK)$ must 
have the form $(F,\KK)$. Then $\mu_{\alpha}(F,\KK) \geq b = \mu_{\alpha}(E_r,\KK)$. If $F \neq E_r$, this contradicts 
the $\alpha^+$-stability of $(E_r,\KK)$. So $F = E_r$ and we have 
$$
0 \ra (\cO(b)^s,0) \ra (G,\KK) \ra (E_r,\KK) \ra 0.
$$
If $s<r$, this contradicts (ii). 

The last assertion is obvious. This completes the proof of (a). 

(b): Suppose that $(E,\KK)$ is $\alpha^-$-stable. Then (ii) is obvious. For the proof of (i) suppose $(F,W) \subset (E_r,\KK)$
contradicts $\alpha^-$-stability. Arguing as in (a) we see that we must have $W = 0$. But then 
$\mu_{\alpha^-}(F,0) > \mu_{\alpha^-}(E,\KK)$, a contradiction.

Suppose conversely, (i) and (ii) hold and $(G,W)$ is a subsystem of $(E,\KK)$ with 
$\mu_{\alpha^-}(G,W) \geq \mu_{\alpha^-}(E,\KK)$. We must have $\mu_{\alpha}(G,W) = \mu_{\alpha}(E,\KK)$
and hence $W =0$. The image of $(G,0)$ in $(\cO(b)^r,0)$ must have the form $(\cO(b)^s,0)$ for some $s \leq r$;
otherwise $\mu_{\alpha}(G,0) < \mu_{\alpha}(E,\KK)$. Hence the intersection of $(G,0)$ and $(E_r,\KK)$ has 
the form $(F,0)$ with $\mu_{\alpha}(F,0) =b = \mu_{\alpha}(E_r,\KK)$. But then 
$\mu_{\alpha^-}(F,0) > \mu_{\alpha^-}(E_r,\KK)$, a contradiction.
\end{proof}

\begin{cor} \label{cor2.2}
{\em (a)} Let $\cM^+$ denote the moduli space of $\alpha^+$-stable coherent systems $(E,\KK)$ which occur in 
an extension {\em (\ref{eq2.1})} and for which $(E, \KK)$ does not admit a coherent subsystem isomorphic to 
$(\cO(b)^{r+1},0)$.
Then $\cM^+$ is isomorphic to a
$Gr(r, C_{12})$-fibration over $G(\alpha^+;n-r,d-br,1) \setminus \cN^+$ 
where $C_{12} =\dim {\rm Ext}^1((E_r,\KK),(\cO(b),0))$ and $\cN^+$ is the moduli space of coherent systems $(E_r, \KK)$ 
which are $\alpha^+$-stable and admit a coherent subsystem $(\cO(b),0)$.

{\em (b)} Let $\cM^-$ denote the moduli space of $\alpha^-$-stable coherent systems $(E,\KK)$ which occur in 
an extension {\em (\ref{e2})} and for which $(E, \KK)$ does not admit a quotient coherent system isomorphic to 
$(\cO(b)^{r+1},0)$.
Then $\cM^-$ is isomorphic to a
$Gr(r, C_{21})$-fibration over $G(\alpha^-;n-r,d-br,1) \setminus \cN^-$ 
where $C_{21} =\dim {\rm Ext}^1((\cO(b),0),(E_r,\KK))$ and $\cN^-$ is the moduli space of coherent systems $(E_r, \KK)$ 
which are $\alpha^-$-stable and admit a quotient coherent system $(\cO(b),0)$.
\end{cor}

\begin{proof} For the proof of (a) just note that in $(\ref{eq2.1})$ the coherent system $(E, \KK)$ admits a subsystem 
$(\cO(b)^{r+1},0)$ if and only if $(E_r,\KK)$ admits a subsystem $(\cO(b),0)$. Similarly, for the proof of
(b)  note that in  $(\ref{e2})$ the coherent system $(E, \KK)$ admits a quotient system 
$(\cO(b)^{r+1},0)$ if and only if $(E_r,\KK)$ admits a quotient system $(\cO(b),0)$.
\end{proof}

For the next proposition we need some notation. Write $a_0 := a, t_0 := t$ and $e_0 := e$ so that 
$$d = na_0 -t_0 \quad \mbox{and} \quad  a_0 = l_0(n-1) + m_0 + t_0
$$ 
with $0 \leq t_0 \leq n-1$, $0 \leq m_0 \leq n-2$ and $0 \leq e_0 \leq l_0$. 
Then define inductively, for $1 \leq r \leq n-2$, 
\begin{itemize} 
\item $s_r$ and $t_r$ by $s_r(n-r) + t_r = t_{r-1} + e_{r-1}$ with $0 \leq t_r \leq n-r-1$,
\item
$a_r := a_{r-1} -s_r$,
\item
$l_r$ and $m_r$ by $a_r = l_r(n-r-1) +m_r +t_r$ with $0 \leq m_r \leq n-r-2$ and
\item
$e_r := e_{r-1} +s_r.$
\end{itemize}   

\begin{lemma} \label{lem3.3}
For $0 \leq r \leq n-2$,\\
{\em (i)} $a_r + e_r = a_0 + e_0$,\\
{\em (ii)} $(n-r)a_r - t_r = d -ra - re$,\\
{\em (iii)} $(n-r)e_r + t_r = ne_0 + t_0 = \alpha_e$,\\ 
{\em (iv)} $0 \leq e_0 \leq e_r \leq l_r - l_0 + e_0 \leq l_r $.\\
In particular, if $1 \leq e \leq l-1$, $\alpha_e$ is a critical value for coherent systems of type $(n-r,d-ra-re,1)$.
\end{lemma}

\begin{proof}
(i) is obvious from the definitions. The proof of (ii), (iii) and (iv) is by induction on $r$,
the case $r=0$ being immediate. So suppose $r \geq 1$ and (ii), (iii) and (iv) are true for $r-1$.
\begin{eqnarray*}
\mbox{(ii)}:\; (n-r)a_r - t_r & = & (n-r)(a_{r-1} - s_r) - (t_{r-1} + e_{r-1} - s_r(n-r))\\
& = & (n-(r-1))a_{r-1} -t_{r-1} - a_{r-1} - e_{r-1}\\
& = &  d - (r-1)a - (r-1)e - a - e \\ 
&  & \;\mbox{(using the induction hypothesis and (i))}\\
& = & d -ra -re.\\
\mbox{(iii)}: \; (n-r)e_r + t_r & = & (n-r)e_{r-1} + (n-r)s_r + t_r\\
& = & (n-r)e_{r-1} + t_{r-1} + e_{r-1}\\
& = & ne_0 + t_0 \quad \mbox{(by induction hypothesis)}.
\end{eqnarray*}
(iv):  By induction hypothesis, $e_r = e_{r-1} + s_r \geq e_0 + s_r \geq e_0$ and on the other hand,
$e_r = e_{r-1} + s_r \leq l_{r-1} - l_0 + e_0 + s_r$. Hence it remains to show that $l_r - l_{r-1} - s_r \geq 0$.
Now we have 
\begin{eqnarray*}
a_{r-1} - s_r = a_r & = &l_r(n-r-1) + m_r + t_r\\
& = & l_r(n-r-1) + m_r + t_{r-1} + e_{r-1} - s_r(n-r).
\end{eqnarray*} 
Hence
$$
l_{r-1}(n-r) + m_{r-1} +t_{r-1} -s_r = l_r(n-r-1) + m_r + t_{r-1} + e_{r-1} - s_r(n-r).
$$ 
Equivalently, 
$$
l_{r-1}(n-r) = l_r(n-r-1) + m_r - m_{r-1} + e_{r-1} -s_r(n-r-1).
$$ 
This gives 
\begin{eqnarray*}
(n-r-1)(l_r - l_{r-1} - s_r) & = & l_{r-1} -m_r + m_{r-1} - e_{r-1} \\
& \geq & l_0 - e_0 -m_r + m_{r-1}\\
&& \mbox{ (by induction hypothesis)}\\
& \geq & -m_r \geq -(n-r-2). 
\end{eqnarray*}
Since $l_r - l_{r-1} - s_r$ is an integer, the left hand side is $\geq 0$ which implies $l_r -l_{r-1} - s_r \geq 0$ and thus completes
the proof.
\end{proof}

Now write for $1 \leq r \leq n-2$
\begin{itemize}
\item $G^r_e := G_{e_r}(n-r,d-ra-re,1)$,
\item $ G_{e-}^r := G_{e_r-1}(n-r,d-ra-re,1)$,
\item $G_e^{r+} := G_{e_r}^+(n-r,d-ra-re,1)$ and
\item $G_e^{r-} := G_{e_r}^-(n-r,d-ra-re,1)$.
\end{itemize}

For $r=n-1$ we write in view of Proposition \ref{prop1.1}
\begin{itemize}
\item $G_e^{n-1} = G_{e-}^{n-1} := G(1,d-(n-1)a - (n-1)e,1)$,
\item $ G_e^{n-1+} = G_e^{n-1-} := \emptyset$.
\end{itemize}

\begin{proposition} \label{prop3.4}
Suppose $1 \leq e \leq l-1$.

{\em (a)}: The variety $G^+_e$ is the disjoint union of the following locally closed subvarieties: 
$$
G_e^+ = \bigsqcup_{r=1}^{n-1} \; V_e^{r+},
$$
where $V_e^{r+}$ is a 
$Gr(r,a-t-(n-1)(e+1)+r)$-bundle over $G_e^r \setminus G_e^{r+}$.  

{\em (b)}: The variety $G^-_e$ is the disjoint union of the following locally closed subvarieties: 
$$
G_e^- = \bigsqcup_{r=1}^{n-1} \;V_e^{r-},
$$
where $V_e^{r-}$ is a 
$Gr(r,r-n +en + t)$-bundle over $G_{e-}^r \setminus G_e^{r-}$. 

{\em (c)}: The base spaces of the bundles defining $V_e^{r+}$ and $V_e^{r-}$ are the same.
\end{proposition}

\begin{proof} (a): Suppose $(E,\KK) \in G_e^+$. By definition there exists a subsystem $(F,W)\subset (E,\KK)$
such that $\mu_{\alpha_e}(F,W) = \mu_{\alpha_e}(E,\KK)$. If $W = \KK$, then $(F,W)$ contradicts $\alpha_e^+$-stability
of $(E,\KK)$. So $W=0$.
Let $\cO(b)$ be a direct factor of $F$ of maximal degree. Replacing $F$ by $\cO(b)$ if necessary, we can
assume without loss of generality that $F = \cO(b)$. 
For equality of 
$\alpha_e$-slopes of $(E,\KK)$ and $(\cO(b),0)$ we must have $b = a+e$. So we get an exact sequence
\begin{equation} \label{eq2.3}
0 \ra (\cO(a+e),0) \ra (E,\KK) \ra (E_1,\KK) \ra 0.
\end{equation}
Suppose that $r \leq n-1$ is the largest integer such that $(E, \KK)$ occurs in an extension
$$
0 \ra (\cO(a+e)^r,0) \ra (E, \KK) \ra (E_r, \KK) \ra 0.
$$ 
By Corollary \ref{cor2.2} (a) the moduli 
space of coherent systems $(E, \KK)$ which occur in this way is isomorphic to a $Gr(r,C_{12})$-bundle over 
$G^r_e \setminus G^{r+}_e$, where $C_{12} = \dim Ext^1((E_r,\KK),(\cO(a+e),0))$. So equation \eqref{e3} implies the assertion.

(b): Suppose $(E,\KK) \in G_e^-$. By definition there exists a subsystem $(F,W) \subset (E,\KK)$
 such that $\mu_{\alpha_e}(F,W) = \mu_{\alpha_e}(E,\KK)$.
If $W = 0$, then $\mu_{\alpha_e^-}(F,W) > \mu_{\alpha_e^-}(E,\KK)$, a contradiction. So $W = \KK$. 
Hence the quotient coherent system is $(G,0)$. Let $\cO(b)$ 
be a direct factor of smallest degree of $G$. Without loss of generality we can assume $G = \cO(b)$. For equality of 
$\alpha_e$-slopes of $(E,\KK)$ and $(\cO(b),0)$ we must have $b = a+e$. So we get an exact sequence
\begin{equation} \label{eq2.2}
0 \ra (E_1,\KK) \ra (E,\KK) \ra (\cO(a+e),0) \ra 0.
\end{equation} 
Suppose that $r \leq n-1$ is the largest integer such that $(E,\KK)$ occurs in an extension
$$
0 \ra (E_r,\KK) \ra (E,\KK) \ra \cO(a+e)^r,0) \ra 0.
$$
By Corollary \ref{cor2.2} (b) the moduli space of coherent systems $(E,\KK)$ which occur in this way is isomorphic 
to a $Gr(r,C_{21})$-bundle over $G^r_{e-1} \setminus G_e^{r-}$, where $C_{21} = \dim Ext^1((\cO(a+e),0),(E_r,\KK))$. 
So equation \eqref{e5} completes the proof. 

(c): This follows at once from equation \eqref{eqn7}.
\end{proof}

We now summarize the results of this section in a theorem.
\begin{theorem} \label{thm3.5}
For $1 \leq e \leq l-1$, the moduli space $G_{e}$ can be obtained from $G_0$ by a series of flips at the critical 
values $\alpha_1, \ldots, \alpha_{e}$. 

The flip at $\alpha_{e'}$ consists of the removal of disjoint locally closed 
subvarieties $V_{e'}^{r-}$ and the insertion of disjoint locally closed subvarieties $V_{e'}^{r+}$ for $1 \leq r \leq n-1$. 
Each of the $V_{e'}^{r+}$ and $V_{e'}^{r-}$ is a Grassmannian bundle over an open subset of a moduli space of coherent systems of 
rank $n-r$. These are described explicitly in Proposition \ref{prop3.4}.
\end{theorem}
\begin{proof}
This follows at once from Proposition \ref{prop3.4}.
\end{proof}

\section{Determination of $G_0$ and $G_{l-1}$}

In order to apply Theorem \ref{thm3.5}, we need to determine $G_0$. In this section we give an inductive description 
of $G_0$ and also of the final moduli space $G_{l-1}$. 
As usual we assume $l \geq 1$, i.e. $a-t \geq n-1$.

\begin{proposition} \label{prop2.1}
If $t=0$, i.e. $d = na$, the moduli space $G_0 = G_0(n,d,1)$ is isomorphic to $Gr(n,a+1)$. 
In particular $G_0 \neq \emptyset$ if and only if $a \geq n-1$.
\end{proposition}

\begin{proof}
If $(E,\KK) \in G_0$, the vector bundle $E$ must be semistable. So $E \simeq \cO(a)^n$.
For $(\cO(a)^n,\KK)$ let $\KK$ be generated by a section 
$\sigma = (\sigma_1, \ldots, \sigma_n) \in H^0(\cO(a))^n$. The coherent system $(\cO(a)^n,\KK)$ is $0^+$-stable if and only if 
$\sigma_1, \ldots, \sigma_n$ are linearly independent. Two such coherent systems are isomorphic if and only if the $n$-tuples 
differ by an element of GL$(n)$. This implies the assertion. 
\end{proof}

When $t>0$, then $G_0 = G_0^+$ and we can use the procedure of Proposition \ref{prop3.4} (a) to determine $G_0$ inductively. 
For this the numbers $s_r, t_r, a_r, l_r, m_r$ and $e_r$ are defined as in section 3. Note that in this case $e_0 = 0$.

\begin{proposition} \label{prop4.2}
For $t > 0$ the variety $G_0$ is the disjoint union of the following locally closed subvarieties: 
$$
G_0 = \bigsqcup_{r=n-t}^{n-1} \; V_0^{r+},
$$
where $V_0^{r+}$ is a 
$Gr(r, a-t+r-n+1)$-bundle over $G_0^r \setminus G_0^{r+}$. 
\end{proposition}

\begin{proof}
Just substitute $e=0$ in the proof of Proposition \ref{prop3.4} (a). We have only to show that $G_0^r = G_0^{r+}$ for 
$1 \leq r \leq n-t-1$.

However $\deg E_r = d-ra = (n-r)a - t$. So if $t \leq n-r-1$, then $G(\alpha;n-r,d-ra,1) = \emptyset$ for $\alpha < t$
which implies the assertion.
\end{proof}

In particular, for $t=1$, we obtain

\begin{cor} \label{cor4.3}
For $t=1$ the moduli space $G_0$ is a $Gr(n-1,a-1)$-bundle over $\PP^{a-1}$.
\end{cor}

\begin{proof}
This follows from the proposition noting that $G_0^{n-1} = G_{e_n-1}(1,d-(n-1)a,1) \simeq \PP^{a-1}$ and 
$G_0^{n-1+}$ is empty.
\end{proof}

In later sections we do not need an explicit description of $G_{l-1}$. However we can easily obtain one.
In fact, by equation \eqref{e1}, $G_{l-1} = G_l^-$ and we can use the procedure of Proposition \ref{prop3.4} (b) to determine $G_{l-1}$
inductively. In this case $e_0 = l$.

\begin{proposition} \label{prop4.3}
The variety $G_{l-1}$ is the disjoint union of the following locally closed subvarieties:
$$
G_{l-1} = \bigsqcup_{r=n-m-1}^{n-1} \; V_{l-1}^{r-},
$$
where $V_{l-1}^{r-}$ is a 
$Gr(r,r-n+ln+t)$-bundle over $G_{l-}^r \setminus G_l^{r-}.$
\end{proposition}

\begin{proof}
Just substitute $e=l$ in the proof of Proposition \ref{prop3.4} (b). We have only to show that $G_{l-}^r = G_l^{r-}$ for 
$1 \leq r \leq n-m-2$.

However by Lemma \ref{lem1},
\begin{eqnarray*}
\deg E_r & = & d - ra -rl \\
& = & (n-1)(ln+t) + mn - r(n-1)l -rm -rt -rl \\
& = & (n-1-r)(ln+t) + m(n-r).
\end{eqnarray*}
So, provided $m \leq n-r-2$, the top limit of $\alpha$ for the existence of $\alpha$-stable coherent 
systems $(E_r,\KK)$ is $ln+t$. This implies $G_{l-}^r = G_l^{r-}$ in this range.
\end{proof}

\section{Hodge polynomials}

Suppose $\KK = \CC$. Then, for any quasiprojective variety $X$, Deligne defined in \cite{d} a mixed Hodge structure on the cohomology groups
$H^k_c(X,\CC)$ with compact support with associated Hodge polynomial $\epsilon(X)(u,v)$. When $X$ is a smooth projective
variety, we have 
$$
\epsilon(X)(u,v) = \sum_{p,q} h^{p,q}(X)u^pv^q
$$
where $h^{p,q}(X)$ are the usual Hodge numbers. In particular, in this case  
$\epsilon(X)(u,u)$ is the usual Poincar\'e polynomial
$P(X)(u)$. We need only the following properties of the Hodge polynomials 
(see \cite{d} and \cite[Theorem 2.2 and Lemma 2.3]{mov}). 

\begin{itemize}
\item If $X$ is a finite disjoint union $X = \sqcup_i X_i$ of locally closed subvarieties $X_i$, then
$$
\epsilon(X) = \sum_{i} \epsilon(X_i).
$$
\item If $Y \ra X$ is an algebraic fibre bundle with fibre $F$ which is locally trivial in the Zariski topology, then 
$$
\epsilon(Y) = \epsilon(X) \cdot \epsilon(F).
$$
\end{itemize}
Moreover we need the Hodge polynomials of the Grassmannians. In fact,
\begin{equation} \label{eqn5.1}
\epsilon(Gr(r,N))(u,v) = \frac{(1 - (uv)^{N-r+1})(1- (uv)^{N-r+2}) \cdots (1- (uv)^{N})}{(1-uv)(1-(uv)^2) \cdots (1-(uv)^r)}.
\end{equation}
Using these properties, we deduce from Propositions \ref{prop3.4}, \ref{prop2.1} and \ref{prop4.2},
\begin{theorem} \label{thm5.1} 
Suppose $\KK = \CC$. Then for any type $(n,d,1)$ with $n \geq 2$ and any integer $e, \; 0 \leq e \leq l-1$, 
the Hodge polynomials $\epsilon(G_e(n,d,1)),
\; \epsilon(G_e^+(n,d,1))$ and $\epsilon(G^-_{e+1}(n,d,1))$ can be explicitly computed. In particular,
$$
h^{p,q}(G_e(n,d,1)) = 0
$$
for $p \neq q$.
\end{theorem}

\begin{proof}
By Propositions \ref{prop3.4} and \ref{prop4.2} and the properties of the Hodge polynomials we have, for $t \geq 1$,
\begin{equation} \label{eqn13}
\epsilon(G_e) = \sum_{r=n-t}^{n-1} \epsilon(V_0^{r+}) + \sum_{e'=1}^e \sum_{r=1}^{n-1} \epsilon(V_{e'}^{r+}) 
- \sum_{e' = 1}^e \sum_{r=1}^{n-1} \epsilon(V_{e'}^{r-}). 
\end{equation}
If $t=0$, the first term on the right hand side of \eqref{eqn13} has to be replaced by $\epsilon(Gr(n,a+1))$ 
according to Proposition \ref{prop2.1}.

The proof proceeds by induction on $n$. For the starting case $n=2$, it follows from Propositions \ref{prop1.1},
\ref{prop3.4}, \ref{prop2.1} and \ref{prop4.2} that $\epsilon(G_e)$ can be expressed in terms of Hodge polynomials of Grassmannians 
and projective spaces. The result follows in this case. 
Now suppose $n \geq 3$ and the theorem is proved for ranks smaller than $n$. The result follows from \eqref{eqn13}
and the description of the varieties $V_{e'}^{r+}$ and $V_{e'}^{r-}$. 
Note that all the fibrations are locally trivial in the Zariski topology, since for $k=1$ there always exist universal
families of coherent systems (see \cite[Appendix]{bgmn}).  

The last assertion follows, since all the ingredients are Grassmannians and projective spaces. 
\end{proof}

In the following sections we will explicitly work out the cases $n=2$ and $n=3$ and also compute $\epsilon(G_0)$ 
for all $n$ when $t \leq 2$.

\section{Coherent systems of rank 2}

Suppose $n = 2$ and as usual $k=1$. Hence $d=2a-t$ and $l = a -t$ with $0 \leq t \leq 1$ and $l \geq 1$.
In particular, for non-emptiness of $G_0$ we require $a-t \geq 1$ and we assume this without further reference.

\begin{proposition} \label{prop5.1}
{\em (a)}: If $t=0$, then $G_0 := G_0(2,d,1)$ is isomorphic to the Grassmannian $Gr(2,a+1)$.\\
{\em (b)}: If $t=1$, then $G_0$ is a $\PP^{a-2}$-bundle over $\PP^{a-1}$.
\end{proposition} 

\begin{proof}
(a) is a special case of Proposition \ref{prop2.1} and (b) is a special case of Corollary \ref{cor4.3}.
\end{proof}

If $a-t=1$, there are no critical values. So suppose $a-t \geq 2$. For the last moduli space 
$G_{a-t-1}$ we have,

\begin{proposition} \label{prop5.2}
The moduli space $G_{a-t-1}$ is isomorphic to $\PP^{2a-t-2}$ and consists of coherent systems $(E,\KK)$ which can 
be expressed as non-trivial extensions
\begin{equation} \label{eq4.1}
0 \ra (\cO,\KK) \ra (E,\KK) \ra (\cO(2a-t),0) \ra 0.
\end{equation}
\end{proposition}

\begin{proof}
This is a special case of Proposition \ref{prop4.3}. In this case $m=0$ and the only allowable value of $r$ is $r=1$.
Now $G_{l-}^1 = G(1,d-a-l,1) \simeq \PP^{d-a-l}$ and $d-a-l = 0$. Moreover $G_l^{1-} = \emptyset$. 
So $G_{a-t-1}$ is isomorphic to $Gr(1,2l+t-1) \simeq \PP^{2l+t-2}$ and $2l+t-2 = 2a -t -2$.
\end{proof}

Our main theorem in this section is a version for genus 0 of results of Thaddeus \cite[section 4]{th}.

\begin{theorem} \label{thm4.2}
Suppose $n=2, \;0\leq t \leq 1$ and $a-t \geq 1$. There exist $a-t$ non-empty moduli spaces $G_0, \ldots, G_{a-t-1}$. 
For $a-t \geq 2$, we have a diagram
\begin{align*}
\begin{small}
\begin{xy}\def\objectstyle{\scriptstyle}
\xymatrix@R=10pt@C=-5pt  { 
& \widetilde{G}_1 \ar[dl] \ar[dr] && \widetilde{G}_2 \ar[dl] \ar[dr] &&\widetilde{G}_3 \ar[dl] \ar[dr]&\;\cdot &\;\cdot &\;\cdot \; 
& \widetilde{G}_{\omega} \ar[dl] \ar[dr] &\\
G_0 && G_1  && G_2 && \;\cdot & \;\cdot &\; \cdot && G_{\omega} .    }\end{xy}
\end{small}
\end{align*}
where $\omega = a-t-1$.

Here $\widetilde{G}_e$ is simultaneously the blow-up of $G_{e-1}$ in $G_e^-$ and the blow-up of $G_{e}$ in $G_{e}^+$.

Moreover, $G_e^+$ is a $\PP^{a-e-t-1}$-bundle over $\PP^{a-e-t}$ and the elements of $G_e^+$ are the coherent systems 
$(E,\KK)$ which can be expressed as non-trivial extensions
$$
0 \ra (\cO(a+e),0) \ra (E,\KK) \ra (\cO(a-t-e),\KK) \ra 0.
$$
The variety $G_e^-$ is a $\PP^{2e+t-2}$-bundle over $\PP^{a-e-t}$ and the elements of $G_e^-$ are the coherent systems 
$(E,\KK)$ which can be expressed as non-trivial extensions
$$
0 \ra (\cO(a-t-e),\KK) \ra (E,\KK) \ra (\cO(a+e),0) \ra 0.
$$
\end{theorem}
\begin{proof}
The assertions about $G_e^+$ and $G_e^-$ are special cases of Proposition \ref{prop3.4} and
the sequences (\ref{eq2.3}) and (\ref{eq2.2}) combined with Proposition \ref{prop1.1}.
Moreover the assumptions A.1 of \cite{bgmn} are satisfied 
for all critical values $\alpha_e = 2e+t$ with $1 \leq e \leq a-t-1$. Hence the diagram exists as claimed. 
\end{proof}

From the explicit description of $G_0$ and $G_{a-t-1}$ 
in Propositions \ref{prop5.1} and \ref{prop5.2} we can compute the Hodge polynomials.

\begin{proposition} \label{prop6.4}
Suppose $\KK = \CC$, $n=2$ and $a-t \geq 1$. Then
$$
\epsilon(G_0)(u,v) = \left\{ \begin{array}{lll} \frac{(1-(uv)^{a})(1-(uv)^{a+1})}{(1-uv)(1-(uv)^2)}& if  & t=0,\\
                                    \frac{(1-(uv)^{a})(1-(uv)^{a-1})}{(1-uv)^2} & if & t=1.
                             \end{array} \right.
$$
$$
\epsilon(G_{a-t-1})(u,v) = \left\{ \begin{array}{lll} \frac{1-(uv)^{2a-1}}{1-uv}& if  & t=0,\\
                                    \frac{1-(uv)^{2a-2}}{1-uv} & if & t=1.
                             \end{array} \right.
$$ 
\end{proposition}

The following lemma gives an inductive description of the Hodge polynomials $\epsilon(G_e)$.

\begin{lemma} \label{lem6.5}
Suppose $\KK = \CC, n=2$ and $1 \leq e \leq a-t-1$. Then
\begin{center}
$\epsilon(G_e)(u,v) - \epsilon(G_{e-1})(u,v) 
= uv \frac{1-(uv)^{a-t-e+1}}{(1-uv)^2}[(uv)^{2e+t-2}- (uv)^{a-t-e-1}].$
\end{center}
\end{lemma}

\begin{proof}
The relevant part of the diagram of Theorem \ref{thm4.2} is
$$
\xymatrix{
& & D \ar[ddll]_{\PP^{a-t-e-1}} \ar[ddrr]^{\PP^{2e+t-2}} \ar@{^{(}->}[d] & & \\
& &       \tilde{G}_e \ar[dl] \ar[dr] & & \\
G_e^- \ar@{^{(}->}[r] & G_{e-1} & & G_e & G_e^+ \ar@{_{(}->}[l]             } 
$$
where $D$ is the exceptional divisor of both blowings-up of $G_e^+$ in $G_e$ and $G_e^-$ in $G_{e-1}$. Hence
$$ 
\epsilon(D) = \epsilon(G_e^+) \cdot \epsilon(\PP^{2e+t-2}) = \epsilon(G_e^-) \cdot \epsilon(\PP^{a-t-e-1}),
$$
\begin{eqnarray*}
\epsilon(\tilde{G}_e) & = & \epsilon(G_{e-1}) + \epsilon(G_e^-)( \epsilon(\PP^{a-t-e-1}) -1)\\
& = &  \epsilon(G_{e}) + \epsilon(G_e^+)( \epsilon(\PP^{2e+t-2}) -1).
\end{eqnarray*}
Now
$$
\epsilon(\PP^b)(u,v) = 1 + uv + \cdots + (uv)^b = \frac{1 - (uv)^{b+1}}{1-uv}.
$$
From the description of $G_e^+$ and $G_e^-$ in Theorem \ref{thm4.2} we get
$$
\epsilon(G_e^+)(u,v) = \frac{1-(uv)^{a-t-e}}{1-uv} \cdot \frac{1-(uv)^{a-t-e+1}}{1-uv},
$$
$$
\epsilon(G_e^-)(u,v) = \frac{1-(uv)^{2e+t-1}}{1-uv} \cdot \frac{1-(uv)^{a-t-e+1}}{1-uv}.
$$
Together this gives the assertion.
\end{proof}
\begin{remark} {\em An alternative proof of Lemma \ref{lem6.5} can be given by observing that 
$\epsilon(G_e) - \epsilon(G_{e-1}) = \epsilon(G_e^+) - \epsilon(G_e^-)$ by equation \eqref{eqn7}
and the first property of Hodge polynomials in section 5}.
\end{remark}

As an immediate consequence we obtain

\begin{cor} \label{cor4.5}
Suppose $\KK = \CC$ and $n=2$. Then
$\epsilon(G_e) = \epsilon(G_{e-1})$ if and only if $3e+2t = a+1$.
\end{cor}

We now have enough information to compute $\epsilon(G_e)$.

\begin{proposition} \label{prop4.5}
Suppose $\KK = \CC$ and $n=2$. Then
$$
\epsilon(G_e)(u,v) = \frac{(1-(uv)^{a-t-e})(1-(uv)^{a-t-e+1})(1-(uv)^{2e+t+1})}{(1-uv)^2(1-(uv)^2)}.
$$
\end{proposition}
\begin{proof}
Write
$$
\epsilon(G_e) = \epsilon(G_0) + \sum_{e'=1}^e (\epsilon(G_{e'}) - \epsilon(G_{e' -1})).
$$ 
Now substitute from Proposition \ref{prop6.4} and Lemma \ref{lem6.5}, use the usual formula 
for summing a geometric series and simplify.
\end{proof}
Substituting $e = a-t-1$, we recover the second formula of Proposition \ref{prop6.4}.

\section{Coherent systems of rank $n \geq 3$}

Suppose $n \geq 3$. For non-emptiness of $G_0$ we require as usual $a-t \geq n-1$ and we 
assume this without further reference. In the same way as in Corollary \ref{cor4.3}, we obtain for $t=2$

\begin{proposition} \label{prop7.1}
If $t=2$, then $G_0 = G_0(n,d,1)$ is a disjoint union of locally closed subvarieties
$$
G_0 = V_0^{n-2+} \sqcup V_0^{n-1+}
$$
where 
\begin{itemize}
\item the variety $V_0^{n-2+}$ is a $Gr(n-2,a-3)$-bundle over $G_0(2,2a-2,1) \, \setminus \, G_1^-(2,2a-2,1)$ 
with $G_0 \simeq Gr(2,a)$ and $G_1^- \simeq \PP^{a-2}$  and
\item the variety $V_0^{n-1+}$ is a $Gr(n-1,a-2)$-bundle over $\PP^{a-2}$.
\end{itemize}
\end{proposition}

Note that $V_0^{n-2+}$ is a dense open subvariety of $G_0$ with complement $V_0^{n-1+}$.

\begin{proof}
By Proposition \ref{prop4.2}, $G_0$ is the disjoint union of locally closed subvarieties 
$V_0^{n-2+}$ and $V_0^{n-1+}$. Here $V_0^{n-2+}$ is a $Gr(n-2,a-3)$-bundle over $G_0^{n-2} \; \setminus \;
 G_0^{n-2+}$ and $V_0^{n-1+}$ is a $Gr(n-1,a-2)$-bundle over $G_0^{n-1} \, \setminus \, G_0^{n-1+}$.
It remains to compute $G_0^{i} \, \setminus \, G_0^{i+}$ for $i=n-1$ and $n-2$.

We have $G_0^{n-1} = G(1,d-(n-1)a,1) \simeq \PP^{d-(n-1)a}$ and $G_0^{n-1+}$ is empty 
according to Proposition \ref{prop1.1}. Moreover $d-(n-1)a = a-2$. For $G_0^{n-2} \, \setminus \, G_0^{n-2+}$ we need 
first to compute $e_{n-2}$.

We claim that $e_{n-2} = 1$. For the proof note that, if $e_{r-1} = 0$ and $t_{r-1} =2$ with $2 < n-r$, then 
$s_r = 0, e_r = 0$ and $t_r = 2$. So the values of these remain constant up to $r = n-3$. Then from 
$2s_{n-2} + t_{n-2} = t_{n-3} + e_{n-3} = 2$ we obtain $s_{n-2} = 1$ and $e_{n-2} = 1$.

So $G_0^{n-2} = G_1(2,d-(n-2)a,1)$ and $d-(n-2)a = 2a-2$. Now by \eqref{eqn7},
$$
G_1(2,2a-2,1) \, \setminus \, G_1^+(2,2a-2,1) = G_0(2,2a-2,1) \, \setminus \, G_1^-(2,2a-2,1).
$$
Finally, $G_0(2,2a-2,1) \simeq Gr(2,a)$ by Proposition \ref{prop2.1}, while $G_1^-(2,2a-2,1)$ is a $Gr(1,1)$-bundle 
over $G(1,a-2,1)$ by Proposition \ref{prop3.4} (b) and $G(1,a-2,1) \simeq \PP^{a-2}$. 
\end{proof}

\begin{cor} \label{cor7.2}
Suppose $\KK = \CC$ and $n \geq 3$. Then the Hodge polynomial of $G_0$ is given by\\
{\em (a)} for $t=0$:
$$
\epsilon(G_0)(u,v) = \frac{(1-(uv)^{a-n+2})(1-(uv)^{a-n+3}) \cdots (1-(uv)^{a+1})}{(1-uv)(1-(uv)^2) \cdots (1-(uv)^n)},
$$
{\em (b)} for $t=1$:
$$
\epsilon(G_0)(u,v) = \frac{(1-(uv)^{a-n+1})(1-(uv)^{a-n+2}) \cdots (1-(uv)^{a})}{(1-uv)^2(1-(uv)^2) \cdots (1 - (uv)^{n-1})},
$$
{\em (c)} for $t=2$:
$$
\epsilon(G_0)(u,v) =  \frac{(1-(uv)^{a-n}) \cdots (1-(uv)^{a-1})(1-(uv)^{n+1})}{(1-uv)(1-(uv)^2)(1-uv)(1-(uv)^2) 
\cdots (1-(uv)^{n-1})}.  
$$
\end{cor}

\begin{proof} (a) and (b) follow from Proposition \ref{prop2.1} and Corollary \ref{cor4.3}.
For (c), from Proposition \ref{prop7.1} we have
$$
\begin{array}{rcl}
\epsilon(G_0)(u,v) & = & \epsilon(V_0^{n-2+}) + \epsilon(V_0^{n-1+})\\
& = & \frac{(1-(uv)^{a-n}) \cdots (1-(uv)^{a-3})}{(1-uv) \cdots (1-(uv)^{n-2})}
[\frac{(1-(uv)^{a-1})(1-(uv)^a)}{(1-uv)(1-(uv)^2)} - \frac{1-(uv)^{a-1}}{1-uv} ]\\
&& + \frac{(1-(uv)^{a-n}) \cdots (1-(uv)^{a-1})}{(1-uv) \cdots (1-(uv)^{n-1}) (1-uv)}.
\end{array}
$$
Simplifying this implies the assertion.
\end{proof}

In order to compute the Hodge polynomial of any moduli space $G_e$, it is now sufficient to work out the change
that takes place at a critical value. In principle this can be done for any $n$. We do it explicitly for $n=3$.

\begin{proposition} \label{prop7.3}
Suppose $\KK = \CC$ and $n = 3$. Then for $1 \leq e \leq l-1$,
\[
\epsilon(G_e) - \epsilon(G_{e-1}) = \frac{1-(uv)^{a-t-2e+1}}{(1-uv)^3(1-(uv)^2)}f(u,v),
\]
where 
\begin{eqnarray*}
f(u,v)&=&(1-uv+(uv)^2)((uv)^{2a-2t-4e-1}-(uv)^{6e+2t-3})\\
&&+
((uv)^{3e+t-2}-(uv)^{a-t-2e-1})(1+(uv)^{a+e+1}).
\end{eqnarray*}
\end{proposition}

\begin{proof}
From equation \eqref{eqn7} and Proposition \ref{prop3.4} we deduce
$$
\begin{array}{rcl}
\epsilon(G_e) &-& \epsilon(G_{e-1})\\
& = & \epsilon(G_e^+) - \epsilon(G_{e}^-)\\
& = & \epsilon(V_e^{1+}) +\epsilon(V_e^{2+}) - \epsilon(V_e^{1-}) - \epsilon(V_e^{2-})\\
& = & (\epsilon(\PP^{a-t-2e-2}) - \epsilon(\PP^{3e+t-3}))(\epsilon(G_e^1) - \epsilon(G_e^{1+}))\\
& + & (\epsilon(Gr(2,a-t-2e)) - \epsilon(Gr(2,3e+t-1)))\epsilon(\PP^{d-2a-2e}).
\end{array}
$$
Recall that $G^1_e:=G_{e_1}(2,d-a-e,1)$ with a similar definition for $G^{1+}_e$. Thus we obtain $\epsilon(G^1_e)$ 
by substituting $a_1$, $t_1$, $e_1$ for $a$, $t$, $e$ in the formula of Proposition 6.8, where $a_1$, $t_1$, $e_1$ 
are defined by
$$2s_1+t_1=t+e\mbox{ with }0\le t_1\le1,\ \ a_1=a-s_1,\ \ e_1=e+s_1$$
(see the definitions preceding Lemma 3.3). Similarly we obtain $\epsilon(G_e^{1+})$ by making the same substitutions 
in the formula for $\epsilon(G_e^+)$ in the proof of Lemma 6.5. Now insert also the Hodge polynomials of 
the Grassmannians and projective spaces and simplify.
\end{proof}

\section{Comments}

\subsection{} A coherent system of type $(n,d,1)$ can be represented by a nonzero homomorphism $\varphi: \cO \ra E$
where $E$ is a vector bundle of rank $n$ and degree $d$, which is determined up to a non-zero scalar multiple.
Such homomorphisms are known as Bradlow pairs. They are special cases of holomorphic triples of rank $(n,1)$ 
and degree $(d,0)$ in the sense of \cite{bgp}. There are appropriate concepts of stability for pairs and for triples 
dependent on a parameter. These parameters for pairs, triples and coherent systems are related by linear
relations. In our case, i.e. $k=1$, the stability conditions then coincide. For $g=0$ the three moduli spaces
are isomorphic.

For $g\geq 2$ and $n=2,3$, results on Hodge polynomials similar to ours have been obtained in \cite{mov} and 
\cite{vm}. For $g=0$, holomorphic triples are discussed in \cite{pp}.

\subsection{} In Corollary \ref{cor7.2} the formulae for $\epsilon(G_0)$ for $t = 0,1,2$ are very simple.
For $t=0$ and $t=1$ this is a direct consequence of the geometric structure of $G_0$. For $t = 2$ the geometric structure 
looks more complicated. One may ask whether there is a simpler description of the geometric structure of $G_0$
which leads naturally to the formulae of Corollary \ref{cor7.2} (c). One may also ask whether there are similarly
simple descriptions of $G_0$ and $\epsilon(G_0)$ for $t \geq 3$.

When $n=2$, the formula for $\epsilon(G_e)$ is also simple (see Proposition 6.8). For $n\ge3$, however, 
it looks as if this formula will be complicated (see Proposition 7.3). It would in any case be good to 
have a geometrical description of $G_e$ which explains the formulae more precisely. Note that as a 
consequence of [5, Theorem 3.2], all moduli spaces $G_e$ are rational varieties.

\subsection{} 
Let $\gamma$ denote the element of the Grothendieck group $K_0(\mbox{Sch}/\mathbb{C})$ of separated 
$\mathbb{C}$-schemes of finite type represented by $G_e(n,d,1)$. Our thanks are due to the referee 
for pointing out that the proof of Theorem 5.1 shows that $\gamma$ belongs to the subring of 
$K_0(\mbox{Sch}/\mathbb{C})$ generated by the affine line. This implies that $\gamma$ is strongly 
polynomial-count in the sense of the appendix by Nicholas M. Katz to \cite{hr}. By \cite[Theorem 6.1.2]{hr}, 
this allows the counting of points of the reduction of $G_e$ over any finite field and hence the 
computation of the zeta function of $G_e$ in terms of the Hodge polynomial. We plan to return to 
this question in a future paper.


\begin{thebibliography}{CAV}

\bibitem{bgp} S. B. Bradlow and O. Garc\'{\i}a-Prada: \emph{Stable triples, equivariant bundles 
and dimensional reduction}. Math. Ann. 304 (1996), 225-252.

\bibitem{bgn} S. B. Bradlow, O. Garc\'{\i}a-Prada, V. Mu\~noz and P. E. Newstead: 
\emph{Coherent systems and Brill-Noether theory}. 
Internat. J. Math. 14 (2003), 683-733.


\bibitem{bgmn} S. B. Bradlow, O. Garc\'{\i}a-Prada, V. Mercat, V. Mu\~noz and P. E. Newstead: 
\emph{On the geometry of moduli spaces of coherent systems on algebraic curves}. 
Internat. J. Math. 18 (2007), 411-453.

\bibitem{d} P. Deligne:
\emph{Theorie de Hodge I,II,III}.
Proc. ICM vol. 1 (1970), 425-430; Publ. Math. IHES 40 (1971), 5-57, ibid. 44 (1974), 5-77.

\bibitem{hr} T. Hausel and F. Rodriguez-Villegas: 
\emph{Mixed Hodge polynomials of character varieties}. 
arXiv:math/0612668v2.

\bibitem{ln} H. Lange and P. E. Newstead: 
\emph{Coherent systems of genus 0}.
Internat. J. Math. 15 (2004), 409-424.


\bibitem{ln1} H. Lange and P. E. Newstead: 
\emph{Coherent systems of genus 0 II: Existence results for $k\geq 3$}.
Internat. J. Math. 18 (2007), 363-393.

\bibitem{vm} V. Mu\~noz: \emph{Hodge polynomials of rank 3 pairs}. arXiv:0706.0593v1.

\bibitem{mov} V. Mu\~noz, D. Ortega and M.-J. V\'azquez-Gallo:
\emph{Hodge polynomials of the moduli spaces of pairs}.
Internat. J. Math. 18 (2007), 695-721.

\bibitem{pp} S. Pasotti and F. Prantil: \emph{Holomorphic triples of genus 0}. To appear in Central European Journ. of Mathem. 

\bibitem{th} M. Thaddeus: 
\emph{Stable pairs, linear systems and the Verlinde formula}.
Invent. Math. 117 (1994), 317-353.

\end{thebibliography}
\end{document}